\documentclass[11pt,draft]{article}
\usepackage{amsmath}
\usepackage{amssymb}
\usepackage{amscd}
\usepackage{xspace}
\usepackage{verbatim}
\newcommand{\mysection}[1]{\section{#1}\setcounter{equation}{0}}
\title{\bf On a new characterisation of Besov spaces with negative exponents}
\author{
{\bf Moshe Marcus}\\
{\small Department of Mathematics,}\\
 {\small Technion, Haifa}
\and
{\bf Laurent V\'eron}\\
{\small Department of Mathematics,}\\
{\small Universit\'e Fran\c{c}ois-Rabelais, Tours}
}

\date{}
\begin{document}%

\maketitle
\noindent{\it Dedicated to Vladimir Maz'ya with high esteem}
\newcommand{\txt}[1]{\;\text{ #1 }\;}
\newcommand{\tbf}{\textbf}
\newcommand{\tit}{\textit}
\newcommand{\tsc}{\textsc}
\newcommand{\trm}{\textrm}
\newcommand{\mbf}{\mathbf}
\newcommand{\mrm}{\mathrm}
\newcommand{\bsym}{\boldsymbol}
\newcommand{\scs}{\scriptstyle}
\newcommand{\sss}{\scriptscriptstyle}
\newcommand{\txts}{\textstyle}
\newcommand{\dsps}{\displaystyle}
\newcommand{\fnz}{\footnotesize}
\newcommand{\scz}{\scriptsize}
\newcommand{\be}{\begin{equation}}
\newcommand{\bel}[1]{\begin{equation}\label{#1}}
\newcommand{\ee}{\end{equation}}
\newcommand{\eqnl}[2]{\begin{equation}\label{#1}{#2}\end{equation}}
\newtheorem{subn}{\name}
\renewcommand{\thesubn}{}
\newcommand{\bsn}[1]{\def\name{#1}\begin{subn}}
\newcommand{\esn}{\end{subn}}
\newtheorem{sub}{\name}[section]
\newcommand{\dn}[1]{\def\name{#1}}   
\newcommand{\bs}{\begin{sub}}
\newcommand{\es}{\end{sub}}
\newcommand{\bsl}[1]{\begin{sub}\label{#1}}
\newcommand{\bth}[1]{\def\name{Theorem}\begin{sub}\label{t:#1}}
\newcommand{\blemma}[1]{\def\name{Lemma}\begin{sub}\label{l:#1}}
\newcommand{\bcor}[1]{\def\name{Corollary}\begin{sub}\label{c:#1}}
\newcommand{\bdef}[1]{\def\name{Definition}\begin{sub}\label{d:#1}}
\newcommand{\bprop}[1]{\def\name{Proposition}\begin{sub}\label{p:#1}}
\newcommand{\R}{\eqref}
\newcommand{\rth}[1]{Theorem~\ref{t:#1}}
\newcommand{\rlemma}[1]{Lemma~\ref{l:#1}}
\newcommand{\rcor}[1]{Corollary~\ref{c:#1}}
\newcommand{\rdef}[1]{Definition~\ref{d:#1}}
\newcommand{\rprop}[1]{Proposition~\ref{p:#1}}
\newcommand{\BA}{\begin{array}}
\newcommand{\EA}{\end{array}}
\newcommand{\BAN}{\renewcommand{\arraystretch}{1.2}
\setlength{\arraycolsep}{2pt}\begin{array}}
\newcommand{\BAV}[2]{\renewcommand{\arraystretch}{#1}
\setlength{\arraycolsep}{#2}\begin{array}}
\newcommand{\BSA}{\begin{subarray}}
\newcommand{\ESA}{\end{subarray}}
\newcommand{\BAL}{\begin{aligned}}
\newcommand{\EAL}{\end{aligned}}
\newcommand{\BALG}{\begin{alignat}}
\newcommand{\EALG}{\end{alignat}}
\newcommand{\BALGN}{\begin{alignat*}}
\newcommand{\EALGN}{\end{alignat*}}
\newcommand{\note}[1]{\textit{#1.}\hspace{2mm}}
\newcommand{\Proof}{\note{Proof}}
\newcommand{\qeda}{\hspace{10mm}\hfill $\square$}
\newcommand{\qed}{\\ ${}$ \hfill $\square$}
\newcommand{\Remark}{\note{Remark}}
\newcommand{\modin}{$\,$\\[-4mm] \indent}
\newcommand{\forevery}{\quad \forall}
\newcommand{\set}[1]{\{#1\}}
\newcommand{\setdef}[2]{\{\,#1:\,#2\,\}}
\newcommand{\setm}[2]{\{\,#1\mid #2\,\}}
\newcommand{\lra}{\longrightarrow}
\newcommand{\lla}{\longleftarrow}
\newcommand{\llra}{\longleftrightarrow}
\newcommand{\Lra}{\Longrightarrow}
\newcommand{\Lla}{\Longleftarrow}
\newcommand{\Llra}{\Longleftrightarrow}
\newcommand{\warrow}{\rightharpoonup}
\newcommand{\paran}[1]{\left (#1 \right )}
\newcommand{\sqbr}[1]{\left [#1 \right ]}
\newcommand{\curlybr}[1]{\left \{#1 \right \}}
\newcommand{\abs}[1]{\left |#1\right |}
\newcommand{\norm}[1]{\left \|#1\right \|}
\newcommand{\paranb}[1]{\big (#1 \big )}
\newcommand{\lsqbrb}[1]{\big [#1 \big ]}
\newcommand{\lcurlybrb}[1]{\big \{#1 \big \}}
\newcommand{\absb}[1]{\big |#1\big |}
\newcommand{\normb}[1]{\big \|#1\big \|}
\newcommand{\paranB}[1]{\Big (#1 \Big )}
\newcommand{\absB}[1]{\Big |#1\Big |}
\newcommand{\normB}[1]{\Big \|#1\Big \|}

\newcommand{\thkl}{\rule[-.5mm]{.3mm}{3mm}}
\newcommand{\thknorm}[1]{\thkl #1 \thkl\,}
\newcommand{\trinorm}[1]{|\!|\!| #1 |\!|\!|\,}
\newcommand{\bang}[1]{\langle #1 \rangle}
\def\angb<#1>{\langle #1 \rangle}
\newcommand{\vstrut}[1]{\rule{0mm}{#1}}
\newcommand{\rec}[1]{\frac{1}{#1}}
\newcommand{\opname}[1]{\mbox{\rm #1}\,}
\newcommand{\supp}{\opname{supp}}
\newcommand{\dist}{\opname{dist}}
\newcommand{\myfrac}[2]{{\displaystyle \frac{#1}{#2} }}
\newcommand{\myint}[2]{{\displaystyle \int_{#1}^{#2}}}
\newcommand{\q}{\quad}
\newcommand{\qq}{\qquad}
\newcommand{\hsp}[1]{\hspace{#1mm}}
\newcommand{\vsp}[1]{\vspace{#1mm}}
\newcommand{\ity}{\infty}
\newcommand{\prt}{\partial}
\newcommand{\sms}{\setminus}
\newcommand{\ems}{\emptyset}
\newcommand{\ti}{\times}
\newcommand{\pr}{^\prime}
\newcommand{\ppr}{^{\prime\prime}}
\newcommand{\tl}{\tilde}
\newcommand{\sbs}{\subset}
\newcommand{\sbeq}{\subseteq}
\newcommand{\nind}{\noindent}
\newcommand{\ind}{\indent}
\newcommand{\ovl}{\overline}
\newcommand{\unl}{\underline}
\newcommand{\nin}{\not\in}
\newcommand{\pfrac}[2]{\genfrac{(}{)}{}{}{#1}{#2}}

\def\ga{\alpha}     \def\gb{\beta}       \def\gg{\gamma}
\def\gc{\chi}       \def\gd{\delta}      \def\ge{\epsilon}
\def\gth{\theta}                         \def\vge{\varepsilon}
\def\vgf{\phi}       \def\vgf{\varphi}    \def\gh{\eta}
\def\gi{\iota}      \def\gk{\kappa}      \def\gl{\lambda}
\def\gm{\mu}        \def\gn{\nu}         \def\gp{\pi}
\def\vgp{\varpi}    \def\gr{\rho}        \def\vgr{\varrho}
\def\gs{\sigma}     \def\vgs{\varsigma}  \def\gt{\tau}
\def\gu{\upsilon}   \def\gv{\vartheta}   \def\gw{\omega}
\def\gx{\xi}        \def\gy{\psi}        \def\gz{\zeta}
\def\Gg{\Gamma}     \def\Gd{\Delta}      \def\vgf{\Phi}
\def\Gth{\Theta}
\def\Gl{\Lambda}    \def\Gs{\Sigma}      \def\Gp{\Pi}
\def\Gw{\Omega}     \def\Gx{\Xi}         \def\Gy{\Psi}

\def\CS{{\mathcal S}}   \def\CM{{\mathcal M}}   \def\CN{{\mathcal N}}
\def\CR{{\mathcal R}}   \def\CO{{\mathcal O}}   \def\CP{{\mathcal P}}
\def\CA{{\mathcal A}}   \def\CB{{\mathcal B}}   \def\CC{{\mathcal C}}
\def\CD{{\mathcal D}}   \def\CE{{\mathcal E}}   \def\CF{{\mathcal F}}
\def\CG{{\mathcal G}}   \def\CH{{\mathcal H}}   \def\CI{{\mathcal I}}
\def\CJ{{\mathcal J}}   \def\CK{{\mathcal K}}   \def\CL{{\mathcal L}}
\def\CT{{\mathcal T}}   \def\CU{{\mathcal U}}   \def\CV{{\mathcal V}}
\def\CZ{{\mathcal Z}}   \def\CX{{\mathcal X}}   \def\CY{{\mathcal Y}}
\def\CW{{\mathcal W}}
\def\BBA {\mathbb A}   \def\BBb {\mathbb B}    \def\BBC {\mathbb C}
\def\BBD {\mathbb D}   \def\BBE {\mathbb E}    \def\BBF {\mathbb F}
\def\BBG {\mathbb G}   \def\BBH {\mathbb H}    \def\BBI {\mathbb I}
\def\BBJ {\mathbb J}   \def\BBK {\mathbb K}    \def\BBL {\mathbb L}
\def\BBM {\mathbb M}   \def\BBN {\mathbb N}    \def\BBO {\mathbb O}
\def\BBP {\mathbb P}   \def\BBR {\mathbb R}    \def\BBS {\mathbb S}
\def\BBT {\mathbb T}   \def\BBU {\mathbb U}    \def\BBV {\mathbb V}
\def\BBW {\mathbb W}   \def\BBX {\mathbb X}    \def\BBY {\mathbb Y}
\def\BBZ {\mathbb Z}

\def\GTA {\mathfrak A}   \def\GTB {\mathfrak B}    \def\GTC {\mathfrak C}
\def\GTD {\mathfrak D}   \def\GTE {\mathfrak E}    \def\GTF {\mathfrak F}
\def\GTG {\mathfrak G}   \def\GTH {\mathfrak H}    \def\GTI {\mathfrak I}
\def\GTJ {\mathfrak J}   \def\GTK {\mathfrak K}    \def\GTL {\mathfrak L}
\def\GTM {\mathfrak M}   \def\GTN {\mathfrak N}    \def\GTO {\mathfrak O}
\def\GTP {\mathfrak P}   \def\GTR {\mathfrak R}    \def\GTS {\mathfrak S}
\def\GTT {\mathfrak T}   \def\GTU {\mathfrak U}    \def\GTV {\mathfrak V}
\def\GTW {\mathfrak W}   \def\GTX {\mathfrak X}    \def\GTY {\mathfrak Y}
\def\GTZ {\mathfrak Z}   \def\GTQ {\mathfrak Q}

\font\Sym= msam10 
\def\SYM#1{\hbox{\Sym #1}}
\newcommand{\bdw}{\prt\Gw\xspace}
\mysection {Introduction}
Let $B$ denote the unit N-ball and $\Sigma=\partial B$. 
If $\mu$ is a distribution on $\Sigma$ we denote by $\BBP (\mu)$ its 
Poisson potential in $B$, that is 
\begin{equation}\label{distripot}\BBP (\mu) (x)=<\mu,P(x,.)>_{\Gs},\quad \forall x\in B,
\end {equation}
where $<\;,\;>_{\Gs}$ denotes the pairing between distributions on $\Gs$ 
and functions in $C^\infty (\Gs)$. In the particular case where $\mu$ is a measure, this can be written 
as follows
\begin{equation}\label{measpot}\BBP (\mu) (x)=\int _{\Sigma}P(x,y)d\mu (y),\quad \forall x\in B.
\end {equation}
In \cite {MV} it is proved that for $q>1$ the Besov space 
$W^{-2/q,q}(\Gs)$ is characterized by an integrability condition on $\BBP 
(\mu)$ with respect to a wheight function involving the distance to 
the boundary, and more precisely that there exists a positive constant 
$C=C(N,q)$ such that for any distribution $\mu$ on $\Gs$ there holds
\begin {equation}\label {old}
C^{-1}{\norm \mu}_{W^{-2/q,q}(\Sigma)}
\leq \left(\int_{B}{\abs {\BBP (\mu)}}^q
(1-\abs x)dx\right)^{1/q}
\leq C{\norm \mu}_{W^{-2/q,q}(\Sigma)}.
\end {equation}
The aim of this article is to prove that for all $1<q<\infty$ any 
 negative Besov spaces $B^{-s,q}(\Gs)$ can be described by an integrability condition on the 
Poisson potential of its elements. More precisely, we prove
\bth {main}Let $s>0$, $q>1$ and $\gm$ be a distribution on $\Gs$. Then 
$$\mu\in B^{-s,q}(\Sigma)\Longleftrightarrow 
\BBP (\mu)\in L^q(B;(1-\abs x)^{sq-1}dx).
$$
Moreover there exists a constant $C>0$ such that for any 
$\gm\in B^{-s,q}(\Sigma)$,
\begin {equation}\label {equivnorm} C^{-1}{\norm \mu}_{B^{-s,q}(\Sigma)}
\leq \left(\int_{B}{\abs {\BBP (\mu)}}^q
(1-\abs x)^{sq-1}dx\right)^{1/q}
\leq C{\norm \mu}_{B^{-s,q}(\Sigma)}.
\end {equation}
\es
The key idea for proving such a result is to use a lifting operator which 
reduces the estimate question to an estimate between Besov spaces with positive 
exponents. In one direction the main technique 
relies on interpolation theory between domain of powers of analytic 
semigroups. In the other direction we use a new representation 
formula for harmonic functions in a ball.
\vspace{5mm}

\noindent{\bf Acknowledgment.}\hspace{2mm}The research of MM was supported by The Israel Science 
Foundation grant No. 174/97.
\mysection {The left-hand side inequality (\ref {equivnorm})}
We recall that for $1\leq p<\infty$, $r\notin \BBN$, $r=k+\eta$ with 
$k\in\BBN$ and $0<\eta<1$, 
$$B^{r,p}(\BBR^d)=\left\{\varphi\in W^{k,p}(\BBR^d):\,
\myint {\BBR^d}{}\myint {\BBR^d}{}\myfrac {{\abs 
{D^\ga\varphi(x)-D^\ga\varphi(y)}}^p}{{\abs {x-y}}^{d+\eta 
p}}dxdy<\infty, \forall \ga\in\BBN^d, \abs \ga=k,\right\}
$$
with norm
$${\norm \varphi}^p_{B^{r,p}}={\norm \varphi}^p_{W^{k,p}}+
\sum_{\abs\ga=k}\myint {\BBR^d}{}\myint {\BBR^d}{}\myfrac {{\abs 
{D^\ga\varphi(x+y)-D^\ga\varphi(x)}}^p}{{\abs {y}}^{d+\eta 
p}}dxdy.
$$
When $r\in \BBN$,
$$\begin {array}{l}B^{r,p}(\BBR^d)=\left\{\varphi\in 
W^{r-1,p}(\BBR^d):^{^{^{^{^{^{}}}}}}\right.\\
\left.\myint {\BBR^d}{}\myint {\BBR^d}{}\myfrac {{\abs 
{D^\ga\varphi(x+2y)+D^\ga\varphi(x)-2D^\ga\varphi(x+y)}}^p}{{\abs 
{y}}^{p+d}}dxdy
<\infty,\forall \ga\in\BBN^d, \abs \ga=r-1,\right\},
\end {array}$$
with norm
$$\begin {array}{l}{\norm \varphi}^p_{B^{r,p}}={\norm \varphi}^q_{W^{k,p}}\\
\qquad\qquad\qquad\qquad +
\displaystyle{\sum_{\abs\ga=r-1}\myint {\BBR^d}{}\myint {\BBR^d}{}\myfrac {{\abs 
{D^\ga\varphi(x+2y)+D^\ga\varphi(x)-2D^\ga\varphi(x+y)}}^p}{{\abs {y}}^{p+d}}dxdy}.
\end {array}$$
The relation of the Besov spaces with integer order of 
differentiation and the classical Sobolev spaces is the following 
\cite {LP}, \cite{Gr}
\begin {equation}\begin {array}{l}\label {inclusion}
B^{r,p}(\BBR^d)\subset W^{r,p}(\BBR^d)\quad \mbox {if } 1\leq p\leq 
2,\\
W^{r,2}(\BBR^d)= B^{r,2}(\BBR^d),\\
W^{r,p}(\BBR^d)\subset B^{r,p}(\BBR^d)\quad \mbox {if }  p\geq 2.
\end {array}\end {equation}
Since for $r\BBN_{*}$ and $1\leq p<\infty$, the space $B^{-r,p}(\BBR^d)$ is 
the space of derivatives of $L^{p}(\BBR^d)$-functions, up to the total 
order $k$, for noninteger $r$, $r=k+\eta$ with $k\in\BBN$ and $0<\eta<1$
$B^{-r,p}(\BBR^d)$ can be defined by using the real interpolation 
method \cite {LP} by
$$\left[W^{-k,p}(\BBR^d),W^{-k-1,p}(\BBR^d)\right]_{\eta,p}=B^{-r,p}(\BBR^d).
$$
The spaces $B^{-r,p}(\BBR^d)$, or $1<p<\infty$ and $r>0$ can also be defined by duality with
$B^{-r,p'}(\BBR^d)$. 
The Sobolev and Besov spaces $W^{k,p}(\Gs)$ and $B^{r,p}(\Gs)$ are 
defined by using local charts from the same spaces in $\BBR^{N-1}$. 
 
\medskip

Now we present the proof of the left-hand side inequality in the case $N\geq 3$.
However, with  minor modifications, the proof 
applies also to the case  $N=2$ (see the remark ). Let $(r,\gs)\in [0,\infty)\times 
S^{N-1}$ (with $S^{N-1}\approx\Gs$) be spherical coordinates in $B$ and put $ t=-\ln r$.
Suppose that $\mu\in B^{-s,q}(S^{N-1})$,
let $u=\BBP(\mu)$ and denote by $ \tl u$  the function $u$ expressed in terms 
of the coordinates $(t,\gs)$. Then 
\begin {equation}\label {equa/u1}
u_{rr}+\myfrac {N-1}{r}u_{r}+\myfrac {1}{r^{2}}\Gd_{\gs}u=0,\quad\mbox 
{in }(0,1)\times S^{N-1},
\end {equation}
and 
\begin {equation}\label {equa/u2}
\tl u_{tt}-(N-2)\tl u_{t}+\Gd_{\gs}\tl u=0,\quad\mbox 
{in }(0,\infty)\times S^{N-1}.
\end {equation}
Then the right inequality in (\ref {equivnorm}) obtains the form
\begin{equation}\label{equivleft}
  \int^\infty_0\int_{S^{N-1}}\abs{\tl 
  u}^q(1-e^{-t})^{sq-1}e^{-Nt}d\gs\,dt\leq 
  C\norm{\mu}^q_{B^{-s,q}(S^{N-1})}.
\end{equation}
Clearly it is sufficient to establish this inequality in the case that $\mu\in \GTM(S^{N-1})$
(or even $\mu\in C^\infty(S^{N-1})$), which is assumed in the sequel. 
We define $k\in\BBN^{*}$ by
\begin{equation}\label{k-def}
    2(k-1)\leq s<2k,
    \end{equation}
    with the restriction $s>0$ if $k=1$. We denote by $\BBb$ 
    the elliptic operator of order 2k
    $$\BBb=\left(\frac{(N-2)^{2}}{4}-\Gd_{\gs}\right)^k
    $$
    and call $f$ the unique solution of 
    $$\gm=\BBb f\quad \mbox {in }S^{N-1}.
    $$
    Then $f\in W^{2k-s,q}(S^{N-1})$ since $\BBb$ is an 
    isomorphism between the spaces $B^{2k-s,q}(S^{N-1})$ and
    $B^{-s,q}(S^{N-1})$.
    Put $v=\BBP(f)$ in $B$, then $v$ satisfies the same 
    equation as $u$ in $(0,1)\times S^{N-1}$. Let $\tl v$
denote this function  in terms of the coordinates $(t,\gs)$. Then 
\begin{equation}\label{equa/v2}
  \begin{cases}
  \tl L{\tl v}:= {\tl v}_{tt}-(N-2){\tl v}_t+\Gd_\gs{\tl v}=0 &\text{in }\BBR_+\ti S^{N-1},\\
  {\tl v}|_{t=0}=f, & \text{in }S^{N-1}.
\end{cases}
\end{equation}
Since the operator $\BBb$ commutes with $\Gd_{\gs}$ and $\partial/\partial t$, 
and this problem has a unique solution which is bounded near $t=\infty$,
it follows that
\begin{equation}\label{equa/v3}
\BBP(\BBb f)=\BBb \tl v.
\end{equation}
Hence,
\begin{equation}\label{equa/v4}
\tl u=\BBP(\mu)=\BBP(\BBb f)
=\BBb \tl v.
\end{equation}
If $v^*:=e^{-t(N-2)/2}\tl v$, then
\begin{equation}\label{equa/v*}\begin{cases}
v^*_{tt}-\frac{(N-2)^2}{4}v^* +\Gd_\gs v^*=0, &\text{in }\BBR_+\ti S^{N-1},\\
v^*(0,\cdot)=f, &\text{in } S^{N-1}.
\end{cases}
\end{equation}
Note that 
$$v^*=e^{tA}(f) \q \text{where} \q 
A=-\paran{\frac{(N-2)^2}{4}I-\Gd_\gs}^{1/2}\Llra A^{2k}=\BBb,$$
where $e^{tA}$ is the semigroup generated by $A$ in $L^q(S^{N-1})$.
By the Lions-Peetre real interpolation method \cite {LP},
$$\left[W^{2k,q}(S^{N-1}),L^q(S^{N-1})\right]_{1-s/2k,q}=B^{2k-s,q}(S^{N-1}).
$$
Since $D(A^{2})=W^{2,q}(S^{N-1})$, $D(A^{2k})=W^{2k,q}(S^{N-1})$.
The semi-group generated by $A$ is analytic as any semi-group 
generated by the square root of a closed operator, therefore
by \cite{Tr} p 96,
\begin{equation}\label{equivleft2}\BA{rcl}
\norm{f}^q_{W^{2k-s,q}}&\sim& \norm{f}_{L^q(S^{N-1})}^q +
\myint{0_{_{}}}{\infty}\paran{t^{(2kqs/2kq)}\norm{A^{2k}v^*}_{L^q(S^{N-1})}}^q\dfrac{dt}{t}\\
&\sim& \norm{f}_{L^q(S^{N-1})}^q +
\myint{0_{_{}}}{1}\paran{t^{s}\norm{A^{2k}v^*}_{L^q(S^{N-1})}}^q\dfrac{dt}{t}\\
&=& \norm{f}_{L^q(S^{N-1})}^q +
\myint{0_{_{}}}{1}\paran{t^{s}e^{-t(N-2)/2}\norm{\BBb\tl v}_{L^q(S^{N-1})}}^q\dfrac{dt}{t}
\EA\end{equation}
where the symbol $\sim$ denotes equivalence of norms. Therefore, by \R{equivleft2},
\begin{equation}\label{equivleft3}\BA{rcl}
    \norm{f}^q_{W^{2k-s,q}(S^{n-1})}&\geq& C \norm{f}_{L^q(S^{N-1})}^q +
C\myint{0_{_{}}}{1}\paran{t^{s}e^{-t(N-2)/2}\norm{\tl u}_{L^q(S^{N-1})}}^q\dfrac{dt}{t}\\
&\geq& C \norm{f}_{L^q(S^{N-1})}^q +
C\myint{0}{1}\norm{\tl u}_{L^q(S^{N-1})}^q e^{-Nt}t^{sq-1}dt.
\EA
\end{equation}
Furthermore,
\begin{equation}\label{equivleft4}\BA{rcl}
\myint{0}{\infty}\norm{\tl u}_{L^q(S^{N-1})}^q(1-e^{t})^{sq-1} e^{-Nt}dt &\leq& C
\myint{0}{1}\norm{\tl u}_{L^q(S^{N-1})}^q (1-e^{-t})^{sq-1} e^{-Nt}dt\\
&\leq& C\myint{0}{1}\norm{\tl u}_{L^q(S^{N-1})}^q e^{-Nt}t^{sq-1}dt.
\EA
\end{equation}
This is a consequence of the inequality
$$\int_{\prt B_r}|u|^q dS\leq(r/\gr)^{N-1} \int_{\prt B_\gr}|u|^q dS,$$
which holds for $0<r<\gr$, for every harmonic function $u$ in $B$. 
By a straightforward computation, this inequality
implies that
$$\int_{|x|<1}|u|^q (1-r) \,dx\leq c(\gg) \int_{\gg<|x|<1}|u|^q (1-r)\,dx,$$
 for every $\gg\in(0,1)$.
 \par In view of  the definition of $f$,
\begin{equation}\label{equivleft5}
\norm{\mu}^q_{B^{-s,q}(S^{n-1})} \sim \norm{f}^q_{W^{2k-s,q}(S^{n-1})}.
\end{equation}
Therefore, the right hand side inequality in \R{equivleft} follows from 
\R{equivleft3},
\R{equivleft4} and \R{equivleft5}.\\[1mm]
\vsp {1}
\mysection {The right-hand side inequality (\ref {equivnorm})}
Suppose that $\mu$ is a distribution on $S^{N-1}$
such that $\BBP(\gm)\in L^q(B;(1-{\abs x}^{sq-1})$. Then we claim that 
$\mu\in B^{-s,q}(S^{N-1})$ and 
\begin{equation}\label{equivright}C^{-1}{\norm \mu}_{B^{-s,q}(\Sigma)}
\leq \left(\int_{B}{\abs {\BBP (\mu)}}^q
(1-\abs x)^{sq-1}dx\right)^{1/q}.
\end {equation}
Because of estimate (\ref{equivleft2}) it is 
suffficient to prove that 
\begin{equation}\label{right1}
 \norm{f}_{L^q(S^{N-1})}\leq C \norm{u}_{L^q(B, (1-r)^{sq-1}\,dx)}.
\end{equation}
With $u=\BBb v$ this relation becomes
\begin{equation}\BA{l}\label{right2}
 \norm{f}_{L^q(S^{N-1})}\leq C \norm{\BBb v}_{L^q(B; 
 (1-r)^{sq-1}\,dx)}\\
 \qquad\;\;\quad\qquad\leq C \left(\myint{0}{1}{\norm 
 v}^q_{W^{2k,q}(S^{N-1})}(1-r)^{sq-1}\,r^{N-1}dr\right)^{1/q}.
\EA \end{equation}
In order to simplify the exposition, we shall first present the case 
where $0<s<2$.
\vsp {1}
\subsection  {\bf The case $0<s<2$}
We take $k=1$. Since the imbedding of 
$B^{2-s,q}(S^{N-1})$ into $L^q(S^{N-1})$ is compact, for any $\vge>0$ 
there is $C_{\vge}>0$ such that
$${\norm \varphi}_{L^q(S^{N-1})}\leq \vge {\norm \varphi}_{B^{2-s,q}(S^{N-1})}+
C_{\vge}{\norm \varphi}_{L^1(S^{N-1})},\quad \forall\,\varphi\in B^{2-s,q}(S^{N-1}).
$$
Therefore the following norm for $B^{2-s,q}(S^{N-1})$ is equivalent to the 
one given in (\ref{equivleft})
\begin{equation}\label{right3}
\norm{f}^q_{B^{2-s,q}}=\norm{f}_{L^1(S^{N-1})}^q +
\myint{0_{_{}}}{1}\paran{t^{s}\norm{A^{2} 
v^{*}}_{L^q(S^{N-1})}}^q\dfrac{dt}{t},
\end{equation}
and estimate (\ref{right2}) will be a consequence of
\begin{equation}\label{right4}
\norm{f}^q_{L^1(S^{N-1})}\leq
\myint{0_{_{}}}{1}\paran{t^{s}\norm{A^{2}
v^{*}}_{L^q(S^{N-1})}}^q\frac{dt}{t}.
\end{equation}
Integrating (\ref {equa/v*}) and using the fact that
\begin{equation}\label{lim}\lim_{t\to\infty}{\norm{v^{*}}}_{L^\infty(S^{N-1})}=
\lim_{t\to\infty}{\norm{v_{t}^{*}}}_{L^\infty(S^{N-1})}=0,
\end{equation}
yields to
$$v_{t}^{*}(t,\gs)=-\myint{t}{\infty}A^{2}v^{*}(s,\gs)ds,\quad 
\forall (t,\gs)\in (0,\infty)\times S^{N-1},
$$
and
\begin{equation}\label{right5}\BA{l}
v^{*}(t,\gs)=\myint{t}{\infty}\myint{s}{\infty}A^{2}v^{*}(\gt,\gs)d\tau ds,\\
\qquad\,\quad\;=\myint{t}{\infty}A^{2}v^{*}(\gt,\gs)(\gt-t)d\gt,
\quad \forall (t,\gs)\in (0,\infty)\times S^{N-1}.
\EA\end{equation} 
Letting $t\to 0$ and integrating over $S^{N-1}$, one obtains
\begin{equation}\label{right6}\BA{l}
\myint{S^{N-1}}{}\abs f d\gs\leq \myint{0}{\infty}\myint{S^{N-1}}{}\abs 
{A^{2}v^{*}}\gt d\gs d\gt \\
\qquad\qquad\;\;\quad\leq 
C(N,s,q,\gd)\left(\myint{0}{\infty}\myint{S^{N-1}}{}{\abs 
{A^{2}v^{*}}}^qe^{\gd \gt}\gt^{sq-1} d\gs d\gt \right)^{1/q}
\EA\end{equation}
for any $\gd>0$ ($\gd$ will be taken smaller that $(N-2)q/2)$ is the 
sequel), where
$$C(N,s,q,\gd)=\left(\abs {S^{N-1}}\myint{0}{\infty}\gt^{(q+1-sq)/(q-1)}e^{-\gd 
\gt/(q-1)}d\gt\right)^{1/q'}.
$$
Notice that the integral is convergent since 
$(q+1-sq)/(q-1)>-1\Longleftrightarrow s<2$.
Going back to $\tl v$
$$\myint{0}{\infty}\myint{S^{N-1}}{}{\abs 
{A^{2}v^{*}}}^qe^{\gd \gt}\gt^{sq-1}d\gs d\gt=
\myint{0}{\infty}\myint{S^{N-1}}{}{\abs {A^{2}\tl v}}^qe^{(\gd-(N-2)q/2) \gt}\gt^{sq-1}d\gs d\gt.$$
Since $u$ is harmonic
$$\myint{S^{N-1}}{}{\abs {\tl u(\gt_{1},.)}}^qd\gs \leq 
\myint{S^{N-1}}{}{\abs {\tl u(\gt_{2},.)}}^qd\gs,\quad \forall 
0<\gt_{2}\leq \gt_{1},
$$
or equivalently,
\begin{equation}\label{right7}\myint{S^{N-1}}{}{\abs {A^{2}\tl v(\gt_{1},.)}}^qd\gs \leq 
\myint{S^{N-1}}{}{\abs {A^{2}\tl v(\gt_{2},.)}}^qd\gs,\quad \forall 
0<\gt_{2}\leq \gt_{1}.
\end{equation}
Applying (\ref{right7}) between $\gt$ and $1/\gt$ for $\gt\geq1$ yields to
\begin{equation}\label{right8}\BA{l}\myint{1}{\infty}\myint{S^{N-1}}{}{\abs 
{A^{2}\tl v }}^qe^{(\gd-(N-2)q/2) \gt}\gt^{sq-1} d\gs d\gt\leq
\myint{0}{1}\myint{S^{N-1}}{}{\abs 
{A^{2}\tl v }}^qe^{(\gd-(N-2)q/2) \gt^{-1}}\gt^{-sq-1} d\gs d\gt
\EA\end{equation}
Moreover there exists $C=C(N,q,\gd)>0$ such that 
$$e^{(\gd-(N-2)q/2) t^{-1}}t^{-sq-1}\leq C
e^{(\gd-(N-2)q/2) t}t^{sq-1}, \quad \forall 0<t\leq 1.$$
Plugging this inequality into \R{right7} and using \R{right6}, one 
derives
\begin{equation}\label{right9}
\myint{S^{N-1}}{}\abs f d\gs\leq C\left(\myint{0}{1}\myint{S^{N-1}}{}{\abs 
{A^{2}v^{*}}}^qe^{\gd \gt}\gt^{sq-1} d\gs d\gt \right)^{1/q}
\end{equation}
for some positive constant $C$, from which \R{right4} follows.
\subsection {the general case}
We assume that $k\geq 1$. Since the imbedding of $B^{2k-s,q}(S^{N-1})$ into 
$L^q(S^{N-1})$ is compact, for any $\vge>0$ 
there is $C_{\vge}>0$ such that
$${\norm \varphi}_{L^q(S^{N-1})}\leq \vge {\norm \varphi}_{B^{2k-s,q}(S^{N-1})}+
C_{\vge}{\norm \varphi}_{L^1(S^{N-1})},\quad \forall\,\varphi\in B^{2k-s,q}(S^{N-1}).
$$
Thus the following norm for $B^{2k-s,q}(S^{N-1})$ is equivalent to the 
one given in (\ref{equivleft})
\begin{equation}\label{right10}
\norm{f}^q_{B^{2k-s,q}}=\norm{f}_{L^1(S^{N-1})}^q +
\myint{0_{_{}}}{1}\paran{t^{s}\norm{A^{2k} 
v^{*}}_{L^q(S^{N-1})}}^q\dfrac{dt}{t},
\end{equation}
and estimate (\ref{right2}) will follow from
\begin{equation}\label{right11}
\norm{f}^q_{L^1(S^{N-1})}\leq
\myint{0_{_{}}}{1}\paran{t^{s}\norm{A^{2k}
v^{*}}_{L^q(S^{N-1})}}^q\frac{dt}{t}.
\end{equation}
From \R{right5}, 
\begin{equation}\label{right12}
    v^{*}(t,\gs)=\myint{t}{\infty}A^{2}v^{*}(\gt,\gs)(\gt-t)d\gt,\quad 
    \forall (t,\gs)\in (0,\infty)\times S^{N-1}.
\end{equation}
Since the operator $A^{2}$ is closed, 
$$
    A^{2}v^{*}(t,\gs)=\myint{t}{\infty}A^{4}v^{*}(\gt,\gs)(\gt-t)d\gt,
$$
and
\begin{equation}\BA{l}\label{right13}
    v^{*}(t,\gs)=\myint{t}{\infty}(t_{1}-t)\myint{t_{1}}{\infty}A^{4}
    v^{*}(t_{2},\gs)(t_{2}-t_{1})dt_{2}dt_{1},\\
    \qquad\quad\;\,=\myint{t}{\infty}\myint{t_{1}}{\infty}(t_{1}-t)(t_{2}-t_{1})A^{4}
    v^{*}(t_{2},\gs)dt_{2}dt_{1},\quad 
    \forall (t,\gs)\in (0,\infty)\times S^{N-1}.
\EA\end{equation}
Iterating this process one gets, for every $ (t,\gs)\in (0,\infty)\times 
S^{N-1}$,
\begin{equation}\label{right14}
v^{*}(t,\gs)=\myint{t}{\infty}\myint{t_{1}}{\infty}\ldots\myint{t_{k-1}}{\infty}
\prod_{j=1}^k(t_{j}-t_{j-1})A^{2k}v^{*}(t_{k},\gs) dt_{k}dt_{k-1}\ldots dt_{1}.
\end{equation}
where we have set $t=t_{0}$ in the product symbol. The following 
representation formula is valid for any $k\in\BBN_{*}$.
\blemma {reduction} For any $(t,\gs)\in (0,\infty)\times S^{N-1}$, 
\begin{equation}\label{right15}
v^{*}(t,\gs)=\myint{t}{\infty}\myfrac{(s-t)^{2k-1}}{(2k-1)!}A^{2k}v^{*}(s,\gs) ds.
\end{equation}
\es
\Proof We proceed by induction. By Fubini's theorem
$$\BA{l}\myint{t}{\infty}\myint{t_{1}}{\infty}(t_{1}-t)(t_{2}-t_{1})A^{4}
 v^{*}(t_{2},\gs)dt_{2}dt_{1}=\myint{t}{\infty}A^{4}
  v^{*}(t_{2},\gs)\myint{t}{t_{2}}(t_{1}-t)(t_{2}-t_{1})dt_{1}dt_{2}\\
   \qquad\quad \qquad\quad \qquad\quad \qquad\quad \qquad\quad \qquad\quad
   =\myint{t}{\infty}\myfrac{(t_{2}-t)^{3}}{6}A^{4}v^{*}(t_{2},\gs)dt_{2}.
\EA$$
Suppose now that for $t>0$, $\ell<k$ and any smooth function $\varphi$ 
defined on $(,\infty)$,
\begin{equation}\label{right16}
\myint{t}{\infty}\myint{t_{1}}{\infty}\ldots\myint{t_{\ell-1}}{\infty}
\prod_{j=1}^\ell(t_{j}-t_{j-1})\varphi(t_{\ell}) dt_{\ell}dt_{\ell-1}\ldots dt_{1}
=\myint{t}{\infty}\myfrac{(t_{\ell}-t)^{2\ell-1}}{(2\ell-1)!}\varphi (t_{\ell}) dt_{\ell}.
\end{equation}
Then
$$\BA{l}\displaystyle{\myint{t}{\infty}\myint{t_{1}}{\infty}\ldots\myint{t_{\ell}}{\infty}
\prod_{j=1}^{\ell+1}(t_{j}-t_{j-1})\varphi(t_{\ell+1}) 
dt_{\ell+1}dt_{\ell}\ldots dt_{1}}\\
\qquad\quad\qquad\qquad\qquad\qquad\quad=
\myint{t}{\infty}\myint{t_{1}}{\infty}\ldots\myint{t_{\ell-1}}{\infty}
\prod_{j=1}^\ell(t_{j}-t_{j-1})\Phi(t_{\ell}) dt_{\ell}dt_{\ell-1}\ldots 
dt_{1},\\
\qquad\quad\qquad\qquad\qquad\qquad\quad=
\myint{t}{\infty}\myfrac{(t_{\ell}-t)^{2\ell-1}}{(2\ell-1)!}\Phi(t_{\ell})dt_{\ell},
\EA$$
with 
$$\Phi(t_{\ell}) =\myint{t_{\ell}}{\infty}(t_{\ell+1}-t_{\ell})\varphi(t_{\ell+1}) dt_{\ell+1}.
$$
But
$$\BA{l}\myint{t}{\infty}\myfrac{(t_{\ell}-t)^{2\ell-1}}{(2\ell-1)!}
\myint{t_{\ell}}{\infty}(t_{\ell+1}-t_{\ell})\varphi(t_{\ell+1}) 
dt_{\ell+1}dt_{\ell}\\
\qquad\qquad\qquad\qquad\qquad\qquad\quad=\myint{t}{\infty}\varphi(t_{\ell+1})
\myint{t}{t_{ell+1}}\myfrac{(t_{\ell}-t)^{2\ell-1}}{(2\ell-1)!}(t_{\ell+1}-t_{\ell})dt_{\ell}dt_{\ell+1}\\
\qquad\qquad\qquad\qquad\qquad\qquad\quad
=\myint{t}{\infty}\varphi(t_{\ell+1})
\myint{0}{t_{ell+1}-t}\myfrac{\gt^{2\ell-1}}{(2\ell-1)!}(t_{\ell+1}-t-\gt)d\gt dt_{\ell+1}\\
\qquad\qquad\qquad\qquad\qquad\qquad\quad
=\myint{t}{\infty}\varphi(t_{\ell+1})\myfrac {(t_{2\ell+1}-\gt)^{2\ell+1}}{(2\ell+1)!}dt_{\ell+1}
\EA$$
as $\myfrac{1}{(2\ell-1)!}(\myfrac {1}{2\ell}-\myfrac 
{1}{2\ell+1})=\myfrac{1}{(2\ell+1)!}$. Taking
$\varphi(t_{\ell+1})=A^{2\ell}v^{*}(t_{\ell+1},\gs)$ implies \R{right15}.
\medskip

\medskip

\noindent {\it End of the proof}. From \R{right14} and Lemma \ref {l:reduction} 
with $t=0$, we get
\begin{equation}\label{right17}\BA{l}
\myint{S^{N-1}}{}\abs f d\gs\leq \myint{0}{\infty}\myint{S^{N-1}}{}\abs 
{A^{2k}v^{*}}\myfrac{\gt^{2k-1}}{(2k-1)!} d\gs d\gt \\
\qquad\qquad\;\;\quad\leq 
C(N,s,k,q,\gd)\left(\myint{0}{\infty}\myint{S^{N-1}}{}{\abs 
{A^{2k}v^{*}}}^qe^{\gd \gt}\gt^{sq-1} d\gs d\gt \right)^{1/q}
\EA\end{equation}
for any $\gd>0$ ($\gd$ will be taken smaller that $(N-2)q/2)$ is the 
sequel), where
$$C(N,s,k,q,\gd)=\left(\abs {S^{N-1}}\myint{0}{\infty}\gt^{(2k-s-1/q')q'}e^{-\gd 
\gt/(q-1)}d\gt\right)^{1/q'}.
$$
Notice that the integral is convergent since 
$(2k-s-1/q')q'>-1\Longleftrightarrow s<2k$. As in the case $s<2$ we return to $\tl 
v$ and $\tl u=A^{2k}\tl u $, use the harmonicity of $u$ in order to 
derive
\begin{equation}\label{right18}\BA{l}\myint{1}{\infty}\myint{S^{N-1}}{}{\abs 
{A^{2k}\tl v }}^qe^{(\gd-(N-2)q/2) \gt}\gt^{sq-1} d\gs d\gt\\
\qquad\qquad\qquad\qquad\qquad\qquad\qquad\qquad\leq
\myint{0}{1}\myint{S^{N-1}}{}{\abs 
{A^{2k}\tl v }}^qe^{(\gd-(N-2)q/2) \gt^{-1}}\gt^{-sq-1} d\gs d\gt
\EA\end{equation}
as in \R{right8} and finally
\begin{equation}\label{right19}\BA{l}
\myint{S^{N-1}}{}\abs f d\gs\leq C\left(\myint {0}{1}\myint{S^{N-1}}{}
{\abs A^{2k}v^{*}}^q\gt^{sq-1}d\gs d\gt \right)^{1/q},\\
\qquad\qquad\qquad \leq C'
\left(\myint {0}{1}\myint{S^{N-1}}{}
{\abs {\tl u}}^q\gt^{sq-1}d\gs d\gt \right)^{1/q},
\EA\end{equation}
which ends the proof of Theorem \ref {t:main}.
\medskip

\noindent {\it Remark}. If $N=2$ the lifting operator is 
$$\BBb=\paran{1-\myfrac {d^{2}}{d\gs^{2}}}^k,
$$
and the proof is similar. moreover, since $\BBb$ is an isomorphism 
between $B^{2k-s,1}(S^{1})$ and $B^{-s,1}(S^{1})$, the 
result of Theorem \ref {t:main} holds also in the case $q=1$.
\mysection {A regularity result for the Green operator}
Put $(1-\abs x)=\gd (x)$. By duality between $L^q(B;\gd^{sq-1}dx)$ and 
$L^{q'}(B;\gd^{sq-1}dx)$, we write
\begin{equation}\label{dual1}
\myint{B}{}\BBP(\mu)\psi \gd^{sq-1}dx=-\myint{B}{}\BBP(\mu)\Gd\zeta dx=-\myint{\Gs}{}\myfrac{\partial 
\zeta}{\partial \gn}d\gm,
\end{equation}
where $\zeta$ is the solution of 
\begin{equation}\label{dual2}\left\{\BA{l}-\Gd\zeta=\gd^{sq-1}\psi\quad \mbox{in }B,\\
    \qquad\zeta=0\qquad \mbox{on }\partial B.\EA\right.\end{equation}
   In \R{dual1}, the boundary term should be written 
    $<\mu,\partial \zeta/\partial \gn>_{\Gs}$ if $\mu$ is a distribution 
    on $\Gs$. Then the adjoint operator $\BBP^{*}$ is defined by
\begin{equation}\label{dual3}
    \BBP^{*}(\psi)=-\myfrac {\partial}{\partial\gn}\BBG(\gd^{sq-1}\psi),
 \end{equation}   
 where $\BBG(\gd^{sq-1}\psi)$ is the Green potential of 
 $\gd^{sq-1}\psi$. Consequently, Theorem \ref {t:main} implies that 
 there exists a constant $C>0$ such that
 
\begin{equation}\label{dual4}
    C^{-1}\norm {\psi}_{L^{q'}(B;\gd^{sq-1}dx)}
    \leq \norm {\myfrac 
    {\partial}{\partial\gn}\BBG(\gd^{sq-1}\psi)}_{B^{s,q'}(\Gs)}
    \leq  C\norm {\psi}_{L^{q'}(B;\gd^{sq-1}dx)}.
 \end{equation} 
 But 
 $$\psi\in L^{q'}(B;\gd^{sq-1}dx)\Longleftrightarrow \gd^{sq-1}\psi\in 
 L^{q'}(B;\gd^{(sq-1)(1-q')}dx).
 $$
 Putting $\varphi=\gd^{sq-1}\psi$ and replacing $q'$ by $p$, implies 
 the following result
 \bth{dual} Let $s>0$ and $1<p<\infty$. Then
 $$\varphi\in L^p(B;\gd^{p(1-s)-1})dx)\Longleftrightarrow
 \myfrac {\partial}{\partial\gn}\BBG(\varphi)\in B^{s,p}(\Gs).
 $$
 Moreover there exists a constant $C>0$ such that for any 
 $\varphi\in L^p(B;\gd^{p(1-s)-1})dx)$
 \begin{equation}\label{dual5}
     C^{-1}{\norm {\varphi}}_{L^p(B;\gd^{p(1-s)-1})dx)}
     \leq {\norm {\myfrac {\partial}{\partial\gn}\BBG(\varphi)}}_{B^{s,p}(\Gs)}
     \leq C{\norm {\varphi}}_{L^p(B;\gd^{p(1-s)-1})dx)}.
  \end{equation}    
 \es

\newpage

\noindent{Moshe Marcus}\\
{\small Department of Mathematics,}\\
 {\small  Israel Institute of Technology-Technion,}\\
 {\small 33000 Haifa, ISRAEL}\\
  {\small Email marcusm@math.technion.ac.il}\\

\noindent {Laurent V\'eron}\\
{\small Laboratoire de Math\'ematiques et Physique Th\'eorique.  CNRS UMR 6083}\\
{\small F\'ed\'eration Denis Poisson }\\
{\small Facult\'e des Sciences, Universit\'e de Tours,}\\
{\small Parc de Grandmont, 37200 Tours, FRANCE}\\
  {\small Email veronl@univ-tours.fr}

\end {document}